\documentclass[12pt]{amsart}
\usepackage{a4wide}
\usepackage{times}
\usepackage{mathtools, amssymb}
\usepackage{graphicx,xspace}
\usepackage{epsfig}
\usepackage{dsfont}
\usepackage[usenames,dvipsnames]{xcolor}
\usepackage{tikz}
\usepackage[T1]{fontenc}
\usepackage[utf8]{inputenc}
\usepackage{array,multirow} % LG : pour faire des beaux tableaux
\usepackage{enumerate}
\usepackage{bm}
\usepackage{cancel}

\usepackage{hyperref,pifont}
\usepackage{todonotes}
\usepackage{dsfont}

\newcommand{\futurework}[1]{}

\usepackage[capitalize]{cleveref}

\theoremstyle{plain}

\newtheorem{lemma}{Lemma}
\newtheorem{theorem}[lemma]{Theorem}

\newtheorem{definition}[lemma]{Definition}

\theoremstyle{remark}

\def\eps{\varepsilon}

\def\R{\mathbb{R}}

\def\tap{\widetilde{\bm a_p}}
\def\tbp{\widetilde{\bm b_p}}
\def\tdp{\widetilde{\bm d_p}}
\def\ttn{\widetilde{\bm t_n}}
\def\ttnj{\widetilde{\bm t_{n,j}}}
\def\ttnp{\widetilde{\bm t_{n,+}}}
\def\ttnm{\widetilde{\bm t_{n,-}}}

\def\pr{\mathbb{P}}
\def\ex{\mathbb{E}}
\DeclareMathOperator{\Var}{Var}
\DeclareMathOperator{\rk}{rk}
\DeclareMathOperator{\des}{des}

\newcommand\restr[2]{{% 
		\left.\kern-\nulldelimiterspace % 
		#1 % 
		\right|_{#2} % 
	}}

\title[CLT for two-sided descents]{On the central limit theorem\\
for the two-sided descent statistics in Coxeter groups}

 \author[V. Féray]{Valentin Féray}
 \address{Institut für Mathematik, Universität Zürich, Winterthurerstr. 190, CH-8057 Zürich, Switzerland}
 \email{valentin.feray@math.uzh.ch}

\keywords{asymptotic normality, Wasserstein distance, Coxeter groups}

\subjclass[2010]{60C05,60F05,05E15}

\begin{document}

\begin{abstract}
  In 2018, Kahle and Stump raised the following problem: identify sequences of finite Coxeter groups $W_n$
  for which the two-sided descent statistics on a uniform random element of $W_n$ is asymptotically normal.
  Recently, Brück and Röttger provided an almost-complete answer, assuming some regularity condition
  on the sequence $W_n$.
  In this note, we provide a shorter proof of their result, which does not require any regularity condition.
  The main new proof ingredient is the use of the second Wasserstein distance on probability distributions,
  based on the work of Mallows (Ann. Math. Statist., 1972).
\end{abstract}

\maketitle

We recall that a sequence of random variables $(X_n)_{n \ge 0}$ is said to be 
{\em asymptotically normal} if $\frac{X_n-\ex[X_n]}{\sqrt{\Var(X_n)}}$
converges in distribution to a standard random variable $Z \sim \mathcal N(0,1)$.
Asymptotic normality of permutation statistics is 
a vast topic in discrete probability, dating back at least to Goncharov \cite{Gon44} and Hoeffding \cite{Hoe51};
we refer also to \cite{Vat96,Ful04,CD17,Ozd19} for more recent works on the descent and two-sided descent statistics.
Recently, there has been some interest into generalizing such asymptotic normality results
to statistics of Coxeter group elements\footnote{For the reader's convenience, we provide an appendix with the necessary definitions regarding Coxeter groups, in particular the notion of {\em descent}.}.
In particular, Kahle and Stump \cite{KS19} have given sufficient and necessary conditions
on a sequence $W_n$ of finite Coxeter groups
so that the number of inversions (resp. of descents) of a uniform random element in $W_n$
is asymptotically normal.
They then asked for a similar characterization for the two sided descent statistics $t$ defined as follows:
for an element $w$ of a Coxeter group $W$, we set $t(w)=\des(w)+\des(w^{-1})$, 
where $\des(w)$ is the number of descents of $w$.
Unlike for inversions and descents, the two sided-descent statistics on a uniform
random element does not decompose as a sum of independent Bernoulli variable,
making the problem more difficult.
For further background on the topic, we refer to \cite{KS19} and \cite{BR19}.
\medskip

The main result of this note is a complete answer to the Kahle--Stump question.
\begin{theorem}
  \label{thm}
  Let $(W_n)_{n \ge 1}$ be a sequence of finite Coxeter groups.
  For each $n$, we let $\bm w_n$ be a uniform random element in $W_n$.  
  Then the following assertions are equivalent:
  \begin{enumerate}[a)]
    \item The sequence $t(\bm w_n)$ is asymptotically normal;
    \item $\Var\big[t(\bm w_n)\big]$ tends to $+\infty$.
  \end{enumerate}
\end{theorem}
This had been previously proved by Brück and Röttger in \cite{BR19}
under a regularity assumption on the sequence $(W_n)_{n \ge 1}$
(the sequence should be {\em well-behaved} in the terminology of \cite{BR19}).
In addition to not requiring any regularity assumption,
the proof that we provide here is shorter.
In particular, we do not need any fourth moment estimates.
\medskip

As in \cite{BR19}, we will take as granted that asymptotic normality holds
when $(W_n)_{n \ge 1}$ is one of the infinite families $A_n$, $B_n$ and $D_n$;
this was proved previously in \cite{Vat96,CD17,Rot18,Ozd19}.
In addition to the fact that $t(w)$ is bounded by $2\rk(W)$,
this is the only specific information we will need on the two-sided descent statistics.
All other arguments are of probabilistic nature.
In particular, we shall use characteristic function analysis,
and Lindeberg type arguments to prove the asymptotic normality (as in \cite{BR19}).
We also introduce a new proof ingredient: 
the second Wasserstein metric for probabilistic distributions.
%The connection between this metric and asymptotic normality
%was suggested in a paper of Mallows \cite{Mal72},
%addressing the asymptotic normality of
%infinite linear combination of independent random variables,
%all of which are asymptotically normal.
%Using the decomposition of Coxeter groups into irreducible components,
%our setting has some similarity with the one studied by Mallows.
\bigskip

We first recall the definition of this Wasserstein metric,
and some useful properties of it, and then proceed to the proof of the main theorem.
\medskip

A (real-valued) random variable $X$ is square integrable if $\ex[X^2]<+\infty$.
A probability distribution (on $\R$) is square integrable 
if a random variable with that probability distribution is.
\begin{definition}[see Lemma 2 in \cite{Mal72}]
Let $\mu$ and $\nu$ be square integrable probability distributions on $\R$. 
Then we define 
\[d_2(\mu,\nu)=\inf_{X \sim \mu, Y \sim \nu} \|X-Y\|_2,\]
where the infimum is taken over all pairs $(X,Y)$ of random variables 
defined on the same probability space and with distributions $\mu$ and $\nu$, respectively.
\end{definition}
As usual in probability theory, we sometimes identify a random variable and its distribution:
namely for random variables $Z$ and $T$ (not necessarily on the same probability space),
we write $d_2(Z,T)= d_2(\pr_Z,\pr_T)$, where $\pr_Z$ and $\pr_T$ are the distributions of $Z$ and $T$.
\medskip

The introduction of the Wasserstein metric (using $L^1$ norm instead of $L^2$ norm, and for general metric space)
is usually attributed to Wasserstein (sometimes also spelled Vasershtein), though
it seems that it appeared in several earlier works \cite{EOM11}.
The $L^2$ case and its relation with asymptotic normality were studied by Mallows \cite{Mal72}.
In particular, he established the following lemmas (Lemmas 1 and 3 in \cite{Mal72}):
\begin{lemma}
  Let $X_n$ and $X$ be square integrable random variables.
  Then $d_2(X_n,X)$ tends to $0$ if and only if
  $X_n \to X$ in distribution and $\ex[X_n^2] \to \ex[X^2]$.
  \label{lem:conv_d2}
\end{lemma}
\begin{lemma}
  Let $k > 0$ be an integer and $Z$ be standard normal random variable.
  If $X_1,\cdots,X_k$ are independent random variables and $(a_j)_{j \le k}$ are real coefficients
  with $\sum_{j \le k} a_j^2 =1$, then
  \[d_2\left( \sum_{j \le k} a_j X_j, Z  \right) \le \sum_{j \le k} a_j^2 \, d_2(X_j, Z).\] 
  \label{lem:add_d2}
\end{lemma}

We can now prove the main result of this note.
\begin{proof}
  [Proof of \cref{thm}]
  The implication a) $\Rightarrow$ b) is immediate: 
  since $t(\bm w_n)$ is integer valued, it cannot tend to a continuous distribution
  without a renormalization factor tending to $+\infty$;
  see \cite[Proposition 6.15]{KS19} for details.
  We focus on b) $\Rightarrow$ a) and assume that $\Var\big[t(\bm w_n)\big]$ tends to $+\infty$.
\medskip

  For each $n \ge 1$, we can decompose the group $W_n$
  as a direct product of irreducible factors $W_n = \prod_{j \le r_n} W_{n,j}$.
  For each $j \le r_n$, we denote by $\bm w_{n,j}$ uniform random elements in $W_{n,j}$
  and by $\bm t_{n,j}=t(\bm w_{n,j})$ the associated two-sided descent statistics.
  Setting $\bm t_n=t(\bm w_{n})$, we have the following decomposition:
  \begin{equation}
  \label{eq:dec_tn}
  \bm t_n\, =\, \sum_{j=1}^{r_n} \bm t_{n,j},
  \end{equation}
  where the $\bm t_{n,j}$ in the right-hand side are independent; see \cite[Lemma 2.2]{BR19}.
  We denote $s^2_{n,j}=\Var\big[t(\bm w_{n,j})\big]$ and $s^2_n=\sum_{j \le r_n} s^2_{n,j} = \Var\big[t(\bm w_n)\big]$.
  Introducing the renormalized random variables 
  \[\ttn:=\frac{\bm t_n - \ex[\bm t_n]}{s_n},\quad
 \ttnj:=\frac{\bm t_{n,j} - \ex[\bm t_{n,j}]}{s_{n,j}},\]
  the decomposition \eqref{eq:dec_tn} writes as
  \[ \ttn = \sum_{j=1}^{r_n} \tfrac{s_{n,j}}{s_n} \ttnj. \]
  Here and in the following, all tilde variables are centered with variance $1$.
  \medskip
  
  We recall that irreducible finite Coxeter groups are of the following types:
  $A_p$ ($p \ge 1$), $B_p$ ($p \ge 2$), $D_p$ ($p \ge 4$), $I_2(m)$ ($m \ge 3$) or
  one of the exceptional types ($H_3$, $H_4$, $E_6$, $E_7$, $E_8$) \cite{Cox35}.
  We write $\bm a_p$, $\bm b_p$ and $\bm d_p$ for uniform random elements in $A_p$, $B_p$
  and $D_p$ respectively.
  As mentioned above, from previous results \cite{Vat96,CD17,Rot18,Ozd19}, 
  we know that the three sequences
  \[\tap:=\frac{\bm a_p - \ex[\bm a_p]}{\sqrt{\Var(\bm a_p)}},\quad
 \tbp:=\frac{\bm b_p - \ex[\bm b_p]}{\sqrt{\Var(\bm b_p)}},\quad
  \tdp:=\frac{\bm d_p - \ex[\bm d_p]}{\sqrt{\Var(\bm d_p)}} \]
  converge in distribution to a standard normal random variable $Z$.
  In addition, their second moment is equal to $1$ for all $p$,
  so we also have convergence of second moments.
  From \cref{lem:conv_d2}, the distributions of $\tap$, $\tbp$ and $\tdp$ converge to that of $Z$
  for the $d_2$ metric.
  \medskip
  
  Fix $\eps >0$ (everything below, including the definitions of large and small components, depends on $\eps$).
   We can find $p_0=p_0(\eps)$ such that for $p \ge p_0$, we have
  \begin{equation}
  \label{eq:d2_le_eps}
   d_2\big(\tap,Z \big) \le \eps,\ d_2\big(\tbp,Z \big) \le \eps,\ d_2\big( \tdp,Z \big) \le \eps.
   \end{equation}
  We now split the irreducible components $(W_{n,j})_{j \le r_n}$ into two groups:
  those of type $A_p$, $B_p$ or $D_p$ for some $p \ge p_0$, which we call {\em large}
  and those of other types 
  to which we will refer to as {\em small}.
  Up to reordering, we can assume that there is an index $q_n=q_n(\eps)$ such that
  large components are exactly those with $j \le q_n$.
\medskip

  We further write $s^2_{n,+}= \sum_{j=1}^{q_n} s^2_{n,j}$ and $s^2_{n,-}= \sum_{j=q_n+1}^{r_n} s^2_{n,j}$.
  We also introduce
  \[\ttnp = \sum_{j=1}^{q_n} \tfrac{s_{n,j}}{s_{n,+}} \ttnj,\quad
  \ttnm = \sum_{j=q_n+1}^{r_n} \tfrac{s_{n,j}}{s_{n,-}} \ttnj,\]
  so that the renormalized two-sided descent statistics decomposes
  into a large component part and a small component part
  \[\ttn = \tfrac{s_{n,+}}{s_{n}}  \ttnp  + \tfrac{s_{n,-}}{s_{n}} \ttnm.\]
  (Summands in the right-hand-side of these equations are independent.)
  \medskip
  
{\em Estimates for the large component part.}
  Using the definition of large components and \cref{eq:d2_le_eps}, we have that 
  $d_2(\ttnj,Z) \le \eps$ for $j \le q_n$. From \cref{lem:add_d2}, this implies
  \[ d_2\big(\ttnp, Z\big) = d_2\left( \sum_{j=1}^{q_n} \tfrac{s_{n,j}}{s_{n,+}} \ttnj, Z \right)
       = \sum_{j=1}^{q_n} \tfrac{s^2_{n,j}}{s^2_{n,+}} d_2 \big( \ttnj, Z \big) 
       \le \eps  \left( \sum_{j=1}^{q_n} \tfrac{s^2_{n,j}}{s^2_{n,+}} \right)
       = \eps. \]
    Using that $u \mapsto \exp(iu)$ is a $1$-Lipschitz function on $\R$,
    we have, for $\zeta$ in $\R$:
    \begin{multline}
          \Big| \ex\big[\exp(i \zeta  \tfrac{s_{n,+}}{s_{n}} \ttnp)\big] 
          - \exp(-  \tfrac{\zeta^2 s^2_{n,+}}{2s^2_{n}})\big] \Big|
          \le \ex\Big[\big|   \exp(i \zeta  \tfrac{s_{n,+}}{s_{n}}\ttnp)  
          -   \exp(i \zeta  \tfrac{s_{n,+}}{s_{n}} Z) \big|\Big] \\
          \le  \tfrac{s_{n,+}}{s_{n}} |\zeta| \, \ex\Big[\big| \ttnp - Z \big|\Big] 
          \le |\zeta| \,  \big\| \ttnp - Z \big\|_2 \le |\zeta| \, \eps,
          \label{eq:est_LargeComp}
    \end{multline}
    where the second to last inequality uses $\tfrac{s_{n,+}}{s_{n}} \le 1$ and Cauchy-Schwartz inequality.
  \medskip
  
{\em Estimates for the small component part.}
  Here, we will use classical characteristic function estimates, as used in Lindeberg central limit theorem
   (see, {\em e.g.}, \cite[Theorem 27.2]{Bil86}).
  By definition, small components are of some exceptional type, 
  of type $I_2(m)$ or of type $A_p$, $B_p$ or $D_p$ for $p<p_0$.
  Their rank is therefore at most $\max(8,p_0)$
  ($I_2(m)$ has rank $2$, the largest exceptional group $E_8$ has rank $8$
  and $A_p$, $B_p$ or $D_p$ have rank $p$).
  But the two sided-descent statistics on any Coxeter group $W$ cannot exceed $2 \rk(W)$.
  We conclude that there is a uniform bound $K=K(\eps)=2\max(8,p_0)$
  on all the $\bm t_{n,j}$ corresponding to small components ($j >q_n$).
  In particular, for $j >q_n$, we have $s_{n,j} \le K$.

  Fix $\zeta$ in $\R$.
  Using, the definition of $\ttnm$, we have
  \begin{equation}
         \ex\Big[\exp\big(i \zeta \tfrac{s_{n,-}}{s_{n}} \ttnm\big)\Big] =
         \prod_{j=q_n+1}^{r_n}   \ex\Big[\exp\big(i  \tfrac{\zeta}{s_{n}} (\bm t_{n,j} - \ex[\bm t_{n,j}]) \big)\Big]\\     
         \label{eq:Tech3}
    \end{equation}  
    We assumed $\lim s_n=+\infty$ and argued above that $s_{n,j}$ is uniformly bounded for $j>q_n$.
   Thus, for $n$ sufficiently large and $j>q_n$, we have $\tfrac{\zeta^2}{s_{n}^2} s_{n,j}^2 \le 1$.
   This implies (see \cite[eqs. (27.11) and (27.15)]{Bil86} that, for $j > q_n$, we have
   \begin{equation}
     \Bigg| \ex\Big[\exp\big(i  \tfrac{\zeta}{s_{n}} (\bm t_{n,j} - \ex[\bm t_{n,j}]) \big)\Big] 
     - \exp\big(-\tfrac{\zeta^2s_{n,j}^2}{2s_{n}^2}\big) \Bigg|
  \le \tfrac{|\zeta|^3}{s_n^3} \ex\Big[ |\bm t_{n,j} - \ex[\bm t_{n,j}]|^3 \Big] + \tfrac{|\zeta|^4}{s_n^4}{s_{n,j}^4} 
  \label{eq:tech_Fourier}
\end{equation}
  Since $\bm t_{n,j}$ is bounded by $K$, we have
  \[\ex\Big[ |\bm t_{n,j} - \ex[\bm t_{n,j}]|^3 \Big] \le K\, \ex\Big[ |\bm t_{n,j} - \ex[\bm t_{n,j}]|^2 \Big]= K s_{n,j}^2.\]
  Using also $s_{n,j}^4 \le K^2 s_{n,j}^2$ and taking $n$ large enough so that $|\zeta| \le s_n$,
  we can simplify the upper bound in \eqref{eq:tech_Fourier} to
  \begin{equation}
    \Bigg| \ex\Big[\exp\big(i  \tfrac{\zeta}{s_{n}} (\bm t_{n,j} - \ex[\bm t_{n,j}]) \big)\Big]
    - \exp\big(-\tfrac{\zeta^2s_{n,j}^2}{2 s_{n}^2} \big) \Bigg|
       \le 2 K^2 |\zeta|^3 \tfrac{s_{n,j}^2}{s_n^3}.  
       \label{eq:tech2}
  \end{equation}
   We now use the following basic inequality: if $(a_i)_{i \le t}$ and $(b_i)_{i \le t}$
   are collections of numbers of absolute values at most one, then
   $\big|\prod_{i \le t} a_i - \prod_{i \le t} b_i\big| \le \sum_{i\le t} |a_i-b_i|$
   (see, {\em e.g.}, \cite[eq. (27.3)]{Bil86}).
   Therefore, \eqref{eq:tech2} implies
   \begin{multline*} \Bigg| \prod_{j=q_n+1}^{r_n}\ex\Big[\exp\big(i  \tfrac{\zeta}{s_{n}} (\bm t_{n,j} - \ex[\bm t_{n,j}]) \big)\Big] 
      - \prod_{j=q_n+1}^{r_n} \exp\big(-\tfrac{\zeta^2s_{n,j}^2}{2 s_{n}^2} \big) \Bigg|\\
      \le \sum_{j=q_n+1}^{r_n} 2 K^2 |\zeta|^3 \tfrac{s_{n,j}^2}{s_n^3} = 2 K^2 |\zeta|^3 \frac{s_{n,-}^2}{s_n^3}
      \le \frac{2 K^2 |\zeta|^3}{s_n}. 
    \end{multline*}
  The first term is the left-hand side is exactly $\ex\Big[\exp\big(i \zeta \tfrac{s_{n,-}}{s_{n}} \ttnm\big)\Big]$;
  see \eqref{eq:Tech3}.
  Since $s_n$ tends to $+\infty$, the upper bound in the last display tends to $0$.
  Therefore for $n$ large enough, we have
  \begin{equation}
    \Bigg| \ex\Big[\exp\big(i \zeta \tfrac{s_{n,-}}{s_{n}} \ttnm\big)\Big] 
    - \exp\big(-\tfrac{\zeta^2 s_{n,-}^2}{2 s_{n}^2} \big)\Bigg| \le \eps.
    \label{eq:est_SmallComp}
 \end{equation}

{\em Conclusion of the proof.}
We recall that $s_{n}^2=s_{n,+}^2+s_{n,-}^2$ and $\ttn = \tfrac{s_{n,+}}{s_{n}}  \ttnp  + \tfrac{s_{n,-}}{s_{n}} \ttnm$.
Using again that $|a_1 a_2 -b_1 b_2| \le |a_1-b_1| +|a_2-b_2|$ for numbers of absolute values at most $1$,
\cref{eq:est_SmallComp,eq:est_LargeComp} imply that, for $n$ large enough,
\[ \Bigg| \ex\Big[\exp\big(i \zeta \ttn\big) \Big] - \exp \big( - \tfrac{\zeta^2}{2} \big) \Big] \Bigg| \le (|\zeta|+1)\,\eps.\]
Since this holds for any $\eps$ and any $\zeta$ in $\R$
(with a threshold value for $n$ depending on $\eps$ and $\zeta$),
we have proved that the characteristic function of $\ttn$ converges pointwise 
towards $\exp \big( - \tfrac{\zeta^2}{2} \big)$,
which is the characteristic function of a Gaussian random variable.
By Lévy's continuity theorem, this concludes our proof.
\end{proof}
\medskip

{\em Technical comment:} a naive characteristic function estimates for the large component
part would lead to an upper bound in \eqref{eq:est_LargeComp} depending on the number $q_n$
of large components. Since we have no control on this number, we would have not been able to conclude.
Using the second Wasserstein distance avoids this problem.

\section*{Appendix: Coxeter groups}

Coxeter groups have been introduced by Coxeter in the '30s \cite{Cox34,Cox35}.
They are now standard objects in combinatorial geometry;
we give here a short introduction to the topic to make this note self-contained.
Classical references are \cite{Hum92,Bou02,BB05}.

A Coxeter matrix $M=(m_{ij})_{i,j \in S}$ indexed by some set $S$
is a symmetric matrix with entries in $\{1,2,3,\cdots\} \cup \{+\infty\}$
such that $m_{ij}=1$ if and only if $i=j$.
A group $W$ is a {\em Coxeter group} if one can find a set $S$ of generators
and a Coxeter matrix $M$ indexed by $S$ such that $W$ admits the presentation
\[W \, \simeq \Big\langle  s \in S \, \big|\, (st)^{m_{st}}=1,\,  s,t \in S \Big\rangle .\]
The pair $(W,S)$ is then called a {\em Coxeter system}.
When we consider a Coxeter group $W$, we often also consider a fixed set $S$
such that $(W,S)$ is a Coxeter system.
The rank of a Coxeter group (or rather of a Coxeter system) is the size of $S$.
Apart from this combinatorial definition, {\em finite} Coxeter groups
can also be characterized geometrically:
they are finite subgroups of general linear groups generated by reflections.

The direct product of two Coxeter groups is a Coxeter group.
A Coxeter group (or rather a Coxeter system)
is irreducible if it cannot be written as a direct product 
of two smaller Coxeter groups.
Trivially, any finite Coxeter group is a direct product of irreducible factors.
Finite irreducible Coxeter groups have been classified by Coxeter in 1935:
\begin{itemize}
  \item there are three infinite families of increasing rank, 
    commonly denoted $A_n$, $B_n$ and $D_n$. 
    the group $A_n$ is the symmetric group on $n+1$ elements,
    $B_n$ is the group of permutations of $n$ elements with 2 colors,
    and $D_n$ is a index 2 subgroup of $B_n$.
  \item there is one infinite family $I_2(m)$ of groups all of rank $2$,
    called {\em dihedral groups}. These are the groups of symmetry of regular polygons.
  \item Finally, there are 6 {\em exceptional groups}, commonly denotes $E_6$, $E_7$, $E_8$,
    $F_4$, $H_3$ and $H_4$ (the index is always the rank or the group).
\end{itemize}
This classification and previous results for the infinite families
are crucial in this note.

We end this appendix by defining the notion of {\em descent} in a Coxeter group
studied in this note.
This generalizes the notion of descents in permutations, corresponding 
to Coxeter groups of type $A_n$.
For an element $w$ in a Coxeter group $W$, we write $\ell(w)$ for the minimal number 
of factors needed to write $w$ as a product of elements of $S$.
Then, by definition, a generator $s$ in $S$ is a descent of $w$ if $\ell(ws)<\ell(w)$.

\section*{Acknowledgements}
I would like to thank Christian Stump for introducing me to the subject of central limit theorems
in Coxeter groups in the "Einstein Workshop on Algebraic Combinatorics" in Berlin, November 2018.
I learned about Mallow's paper \cite{Mal72} through the {\em Encyclopedia of Mathematics}
page on Wasserstein metric \cite{EOM11}.

\end{document}